\documentclass[11pt,british]{amsart}

\usepackage[T1]{fontenc}

\usepackage{babel,eucal,url,amssymb,enumerate,amscd,graphics}

\theoremstyle{plain}
\newtheorem{lemma}{Lemma}[section]

\newtheorem{proposition}[lemma]{Proposition}
\newtheorem{corollary}[lemma]{Corollary}
\newtheorem{theorem}[lemma]{Theorem}
\newtheorem{remark}[lemma]{Remark}

\newcommand{\Gtwo}{\ifmmode{{\rm G}_2}\else{${\rm G}_2$}\fi}

 \newcommand{\cyclic}{\mathop{\kern0.9ex{{+}\kern-2.2ex\raise-.28ex\hbox{\Large\hbox
 {$\circlearrowright$}}}}}

\newcommand{\End}{\mathop{\mbox{\rm End}}}
\newcommand{\Aut}{\mathop{\mbox{\rm Aut}}}

\def\sideremark#1{\ifvmode\leavevmode\fi\vadjust{\vbox to0pt{\vss
 \hbox to 0pt{\hskip\hsize\hskip1em
 \vbox{\hsize2.5cm\tiny\raggedright\pretolerance10000
 \noindent #1\hfill}\hss}\vbox to8pt{\vfil}\vss}}}%

\newfont{\eusm}{eusm10 scaled \magstep1}
\newfont{\eusmiii}{eusm10 scaled \magstep3}

\newcommand{\comp}{\makebox[7pt]{\raisebox{1.5pt}{\tiny $\circ$}}}
\newcommand{\RR}{\mbox{{\sl I}}\!\mbox{{\sl R}}}

\title[Jacobi osculating rank and isotropic geodesics]{Jacobi osculating rank and isotropic geodesics on naturally reductive 3-manifolds}

\author{J.~C.~Gonz{\'a}lez-D{\'a}vila}
\address{Department of Fundamental Mathematic\\
  University of La Laguna\\ 38200 La Laguna, Tenerife, Spain}
\email{jcgonza@ull.es}

\begin{document}
\maketitle

\begin{abstract}{\indent} We study the Jacobi osculating rank of
geodesics on naturally reductive homogeneous manifolds and we
apply this theory to the $3$-dimensional case. Here, each
non-symmetric, simply connected naturally reductive $3$-manifold
can be given as a principal bundle $M^{3}(\kappa,\tau)$ over a
surface of constant curvature $\kappa,$ such that the curvature of
its horizontal distribution is a constant $\tau>0,$ with
$\tau^{2}\neq \kappa.$ Then, we prove that the Jacobi osculating
rank of every geodesic of $M^{3}(\kappa,\tau)$ is two except for
the Hopf fibers, where it is zero. Moreover, we determine all
isotropic geodesics and the isotropic tangent conjugate locus.

\vspace{4mm}

 \vspace{2mm}

\end{abstract}
\begin{figure}[b]  \vspace{-5mm}

\end{figure}

\tableofcontents



\section{Introduction}\indent

A Jacobi field $V$ on a homogeneous Riemannian manifold $(M,g)$
which is the restriction of a Killing vector field along a
geodesic is called {\em isotropic} \cite{Z}. From the homogeneity
of $(M,g),$ it means that $V$ is the restriction of a fundamental
vector field of some element $X$ in the Lie algebra of the
isometry group $I(M,g)$ of $(M,g).$ If moreover $V$ vanishes at a
point of the geodesic, then $X$ belongs to the Lie algebra of the
isotropy subgroup at this point. This particular situation was
what originally motivated the term {\em isotropic} (see \cite{Ch}
and \cite{Ch1}).

Two points $p,q\in M$ are said to be {\em isotropically conjugate}
if there exists a nonzero isotropic Jacobi field $V$ along a
geodesic passing through $p$ and $q$ such that $V$ vanishes at
these points. Clearly, they are isotropically conjugate along any
geodesic joining $p$ to $q.$ When every Jacobi field vanishing at
$p$ and $q$ is isotropic, we say that they are {\em strictly
isotropic conjugate points.}

On symmetric spaces, the Jacobi equation has simple solutions and
one directly obtains that {\em any pair of conjugate points in a
Riemannian symmetric space are strictly isotropic} (see
\cite{GS}). On \cite[Remark 4.7]{GD} are constructed examples of
$3$-symmetric spaces admitting pairs of conjugate points which are
isotropic but not strictly isotropic. A geodesic starting at a
point $p\in M$ is said to be {\em isotropic} (resp., {\em strictly
isotropic}) if each one of its conjugate points to $p$ is
isotropic (resp., strictly isotropic). Then, all geodesic on a
symmetric space is strictly isotropic. In the case of a naturally
reductive space, W. Ziller in \cite{Z} proposed to examine
conjectures like: {\em A naturally reductive space with the
property that all its geodesics are strictly isotropic is locally
symmetric}. A positive answer to the conjecture for $n\leq 5$ is
given in \cite{GS} and in \cite{GD}, for naturally reductive
compact 3-symmetric spaces.

The parallel translation of the Jacobi operator $R_{\gamma_{u}}:=
R(\gamma'_{u},\cdot)\gamma'_{u}$ along a geodesic $\gamma_{u}$
starting at the origin $o$ of a naturally reductive homogeneous
manifold $(M=G/K,g)$ with $\gamma'_{u}(0) =u,$ for some unit
vector $u\in {\mathfrak m}\cong T_{o}M,$ determines a curve
$R_{u}(t)$ in the space of the self-adjoint operators ${\mathcal
S}({\mathfrak m})$ of ${\mathfrak m},$ which is an orbit of a
one-parameter subgroup of isometries of this space. Then, as it is
shown in Lemma \ref{losc}, such a curve has constant osculating
rank. We refer to such constant as the {\em Jacobi osculating rank
of the geodesic} $\gamma_{u}$ and it will be denoted by ${\rm
rank}_{\rm osc}(u).$ Then, there exist smooth functions
$a_{1},\dots ,a_{r},$ where $r = {\rm rank}_{\rm osc}(u),$ such
that
\[
R_{u}(t) = R_{u}(0) + a_{1}(t)R'_{u}(0) + \dots +
a_{r}(t)R_{u}^{r)}(0),
\]
$R'_{u}(0),\dots , R_{u}^{r)}(0)$ are linearly independent in
${\mathcal S}({\mathfrak m})$ and $r\leq\frac{n(n-1)}{2},$ $n =
\dim M.$ Moreover, in Lemma \ref{pcharac} one obtains that the
Jacobi osculating rank is also the osculating rank of the curve
obtained in ${\mathcal S}({\mathfrak m})$ by parallel translation
of the Jacobi operator for the {\em canonical connection} adapted
to a reductive decomposition of $(M=G/K,g).$

When the Jacobi osculating rank does not depend on the choice of
the geodesic, the naturally reductive homogeneous space is said to
have {\em constant Jacobi osculating rank}. Clearly, symmetric
spaces have Jacobi constant osculating rank zero, or equivalently
all these curves are constant, and for the simply connected case,
in Theorem \ref{sim} the converse is stated. Then the notion of
Jacobi osculating rank can be view as a natural way to `measure'
what a naturally reductive homogeneous space, or more general a
g.o. space (see Remark \ref{go}) moves away from to be locally
symmetric.

Some examples of non-symmetric naturally reductive spaces with
constant Jacobi osculating rank are already known. Concretely, A.
M. Naveira and A. Tarr\'{\i}o in \cite{N-T} and, together with E.
Mac\'{\i}as, in \cite{M-N-T} have proved that the Berger manifold
$V_{1} = Sp(2)/SU(2)$ and the Wilking manifold $V_{3}=
(SO(3)\times SU(3))/U^{\bullet}(2),$ endowed this last one with a
particular bi-invariant metric, have constant Jacobi osculating
rank two and, in the context of g.o. spaces, T. Arias-Marco and A.
M. Naveira in \cite{AM-N} have shown that the Jacobi osculating
rank of the six-dimensional Kaplan's example is constant equals
four. (See \cite{GN}, for a brief survey on isotropic Jacobi
fields and Jacobi osculating rank and for further references.) For
a general naturally reductive homogeneous space, the constancy of
the Jacobi osculating rank is not necessarily satisfied. In fact,
as a consequence from Proposition \ref{nonconst}, one obtains that
non-locally symmetric naturally reductive spaces of dimension
$n\leq 5,$ generalized Heisenberg groups, Berger spheres or
$\varphi$-symmetric spaces are some examples of non-constant
Jacobi osculating rank.

In this article, we shall focus our attention on naturally
reductive homogeneous spaces of dimension $3.$ The classification
for the simply connected case is well known \cite{TV} (see also
\cite{Dan}, \cite{K1}). They are the symmetric spaces $\RR^{3},$ $
S^{3}(c),$ $H^{3}(-c),$ $S^{2}(c)\times \RR$ and $H^{2}(-c)\times
\RR,$ where $c>0,$ and unimodular Lie groups equipped with a left
invariant metric such that the dimension of their isometry groups
is four. Each one of these spaces fibers as an one-dimensional
principal fiber bundle over a complete simply connected surface of
constant curvature $\kappa$ and the horizontal distribution of
this fibration is the kernel of a connection form with constant
curvature $\tau,$ which can be taken positive (see \cite{GGV},
\cite{Thu}). The fibers are geodesics and there exists a
one-parameter family of translations along the fibers, generated
by a unit Killing vector field $\xi,$ the {\em Hopf vector field}.
Then these manifolds, which will be denoted by
$M^{3}(\kappa,\tau),$ are classified, up to isometry, in terms of
the pair $(\kappa,\tau),$ where $\kappa\neq \tau^{2}$ and
$\tau>0,$ into three types: Berger spheres $S^{3}(\kappa,\tau),$
if $\kappa>0;$ the universal covering
$\widetilde{SL_{2}}(\kappa,\tau)$ of $SL(2,\RR)(\kappa,\tau),$ if
$\kappa <0;$ and the Heisenberg group $H_{3}(\tau),$ if $\kappa
=0.$

Because $\xi$ is a unit Killing vector field, every geodesic
$\gamma_{u}$ on $M^{3}(\kappa,\tau)$ intersects each fiber with a
constant slope angle $\theta= {\rm ang}(\xi_{o},u)\in ]0,\pi[.$
Moreover, according with Lemma \ref{lad}, every pair of geodesics
starting at the origin with same slope angle are related under the
isotropy action. It allows us to prove in Theorem \ref{tosc} that
{\em the Jacobi osculating rank of every geodesic on $M^{3}
(\kappa,\tau)$ is two except for the Hopf fibers, where it is
zero}. Moreover, for each slope angle $\theta\in ]0,\pi[,$ the
curve $R_{u}(t)$ in ${\mathcal S}({\mathfrak m})$ is a circle
whose radius contracts to the point $R_{\xi_{o}}$ as $\theta$
converges to $0$ or to $\pi.$

S. Engel determined in \cite{Engel} the conjugate radius and the
cut locus of $M^{3}(\kappa,\tau),$ as generalisation of the
results on Berger spheres carried out by Sasaki \cite{Sakai} and
Rakotoniaina \cite{Ra}. Here, we go forward with this research.
First, taking into account that, by using the adapted canonical
connection, the Jacobi equation can be expressed as a differential
equation with constant coefficients, we obtain in Theorem
\ref{pral1} all isotropic conjugate points of
$M^{3}(\kappa,\tau).$ Then its isotropic geodesics can be
explicitly determined. They are one-to-one geodesics without pairs
of conjugate points and with slope angle $\theta$ equals to
$\pi/2$ for $H_{3}(\tau)$ and $\theta\in [\varepsilon, \pi -
\varepsilon],$ where $\varepsilon = \arctan \tau/\sqrt{-k},$ for
$\widetilde{SL_{2}}(\kappa,\tau).$ $S^{3}(\kappa,\tau)$ does not
admit any isotropic geodesic. Finally, the tangent conjugate and
isotropic conjugate locus of $M^{3}(\kappa,\tau)$ are given as
union of surfaces of revolution about the Hopf direction.

\section{Preliminaries}\indent

Let $(M,g)$ be a connected homogeneous {\em Riemannian} manifold.
As it is well-known, $(M,g)$ can be expressed as coset space
$G/K,$ where $G$ is a connected Lie group of isometries acting
transitively and effectively on $M,$ $K$ is the isotropy subgroup
of $G$ at some point $o\in M,$ the {\em origin} of $M,$ and $g$ is
considered as a $G$-invariant Riemannian metric on $G/K.$
Moreover, we can assume that $G/K$ is a {\em reductive homogeneous
space}, i.e., there is an $Ad(K)$-invariant subspace ${\mathfrak
m}$ of the Lie algebra ${\mathfrak g}$ of $G$ such that
${\mathfrak g} = {\mathfrak m} \oplus {\mathfrak k},$ ${\mathfrak
k}$ being the Lie algebra of $K.$ Such quotient representation of
$(M,g)$ is in general not unique. $(M=G/K,g)$ is said to be {\em
naturally reductive}, or more precisely $G$-{\em naturally
reductive}, if there exists a reductive decomposition ${\mathfrak
g} = {\mathfrak m} \oplus {\mathfrak k}$ satisfying
\begin{equation}\label{nred}
<[X,Y]_{\mathfrak m},Z> + <[X,Z]_{\mathfrak m},Y> = 0
\end{equation}
\noindent for all $X,Y,Z\in {\mathfrak m},$ where
$[X,Y]_{\mathfrak m}$ denotes the ${\mathfrak m}$-component of
$[X,Y]$ and $<,>$ is the metric induced by $g$ on ${\mathfrak m},$
by using the canonical identification ${\mathfrak m}\cong T_{o}M.$
If there exists a bi-invariant metric $B$ on ${\mathfrak g}$ whose
restriction to ${\mathfrak m} = {\mathfrak k}^{\bot}$ is the
metric $<,>,$ the (naturally reductive) space $(M=G/K,g)$ is
called {\em normal homogeneous}. Then, for all $X,Y,Z\in
{\mathfrak g},$ we have
\begin{equation}\label{tt2}
B([X,Y],Z) + B([X,Z],Y) = 0.
\end{equation}

  For each $X\in {\mathfrak g},$ the mapping $\psi:\RR \times M \to M,$
$(t,p)\in \RR\times M \mapsto \psi_{t}(p) = (\exp tX)p$ is a
one-parameter group of isometries and consequently, $\psi$ induces
a Killing vector field $X^{*}$ given by
\begin{equation}\label{estrella}
X^{*}_{p} = \frac{\textstyle d}{\textstyle dt}_{\mid t = 0}(\exp
tX)p, \;\;\;\;\; p\in M.
\end{equation}
$X^*$ is called the {\em fundamental vector field} or the {\em
infinitesimal $G$-motion} corresponding to $X$ on $M.$ If $G =
I_{o}(M,g),$ then all (complete) Killing vector field on $M$ is a
fundamental vector field $X^{*},$ for some $X\in {\mathfrak g}.$
For any $a\in G,$ we have
\begin{equation}\label{adjun1}
(Ad_{a}X)^{*}_{ap} = a_{*p}X^{*}_{p},
\end{equation}
where $a_{*p}$ denotes the differential map of $a$ at $p\in M.$

Next, let $\tilde{T}$ denote the torsion tensor and $\tilde{R}$
the corresponding curvature tensor of the {\em canonical
connection} $\tilde{\nabla}$ of $(M,g)$ adapted to the reductive
decomposition ${\mathfrak g} = {\mathfrak m} \oplus {\mathfrak k}$
\cite[I, p.110]{KN} defined by the sign convention $\tilde{R}(X,Y)
= \tilde{\nabla}_{[X,Y]} -
[\tilde{\nabla}_{X},\tilde{\nabla}_{Y}]$ and $\tilde{T}(X,Y) =
\tilde{\nabla}_{X}Y - \tilde{\nabla}_{Y}X - [X,Y],$ for all
$X,Y\in {\mathfrak X}(M),$ the Lie algebra of smooth vector fields
on $M.$ Then, these tensors are given by
\begin{equation}\label{torsion}
\tilde{T}_{o}(X,Y) = -[X,Y]_{\mathfrak m} \;\;\;, \;\;\;
\tilde{R}_{o}(X,Y) = {\rm ad}_{[X,Y]_{\mathfrak k}}
\end{equation}
and we have $\tilde{\nabla}g = \tilde{\nabla}\tilde{T} =
\tilde{\nabla}\tilde{R} = 0.$ On naturally reductive homogeneous
manifolds $(M=G/K,g),$ the tensor field $S = \nabla -
\tilde{\nabla},$ where $\nabla$ denotes the Levi Civita connection
of $(M,g),$ is a {\em homogeneous structure} \cite{TV} satisfying
$S_{X}Y =-S_{Y}X = -\frac{1}{2}\tilde{T}(X,Y),$ for all $X, Y\in
{\mathfrak X}(M),$ and we get
\begin{equation}\label{RR}
\tilde{R}_{XY} = R_{XY} + [S_{X},S_{Y}] - 2S_{S_{X}Y}.
\end{equation}
Then $\nabla$ and $\tilde{\nabla}$ have the same geodesics and,
consequently, the same Jacobi fields (see \cite{Z}). Such
geodesics are orbits of one-parameter subgroups of $G$ of type
$\exp tu$ where $u\in {\mathfrak m}.$ In what follows, we shall
denote by $\gamma_{u}$ the unit-speed geodesic starting at the
origin $o$ with $\gamma_{u}'(0) = u,$ $\|u\|=1.$ Then
$\gamma_{u}(t) := (\exp tu)o$ and the Jacobi equation for $\nabla$
coincides with the Jacobi equation for $\tilde{\nabla},$
\[
\frac{\tilde{\nabla}^{2}V}{dt^{2}} -
\tilde{T}_{\gamma}\frac{\tilde{\nabla}V}{dt} + \tilde{R}_{\gamma}V
= 0,
\]
where $\tilde{R}_{\gamma} = \tilde{R}(\gamma',\cdot)\gamma'$ and
$\tilde{T}_{\gamma} = \tilde{T}(\gamma',\cdot).$ Taking into
account that $\tilde{\nabla}\tilde{T} = \tilde{\nabla}\tilde{R} =
0$ and the parallel translation with respect to $\tilde{\nabla}$
of tangent vectors at the origin along $\gamma_{u}$ coincides with
the differential of ${\exp}tu \in G$ acting on $M,$ it follows
that any Jacobi field $V$ along $\gamma(t)$ can be expressed as
$V(t) = (\exp tu)_{*o}X(t)$ where $X(t)$ is solution of the
differential equation
\begin{equation}\label{Jam1}
X''(t) - \tilde{T}_{u}X'(t) + \tilde{R}_{u}X(t) = 0
\end{equation}
in the vector space ${\mathfrak m},$ being $\tilde{T}_{u}X =
\tilde{T}(u,X) = -[u,X]_{\mathfrak m}$ and $\tilde{R}_{u}X =
\tilde{R}(u,X)u = [[u,X]_{\mathfrak k},u]$ (see \cite{GS},
\cite{Z} for more details).

A Jacobi field $V$ along $\gamma_{u}$ with $V(0) = 0$ is said to
be $G$-{\em isotropic}, or simply {\em isotropic} if
$G=I_{o}(M,g),$ if and only if there exists $A\in {\mathfrak k}$
such that $V = A^{*}\comp \gamma,$ or equivalently, if there
exists an $A\in {\mathfrak k}$ such that (see \cite{GS})
\begin{equation}\label{ini}
V'(0) = [A,u].
\end{equation}

\section{Jacobi osculating rank}\indent

Let $(M = G/K,g)$ be a connected naturally reductive homogeneous
Riemannian manifold with adapted decomposition ${\mathfrak g} =
{\mathfrak m}\oplus {\mathfrak k}.$ For each $t\in \RR,$ denote by
$R_{u}(t)$ the $(1,1)$-tensor on ${\mathfrak m}$ obtained by the
parallel translation of the Jacobi operator $R_{\gamma_{u}}$ along
the geodesic $\gamma_{u},$ i.e.,
\[
R_{u}(t) = \tau_{u}^{-t}\comp R_{\gamma_{u}}\comp \tau^{t}_{u},
\]
where $\tau^{t}_{u}:T_{o}M\cong {\mathfrak m}\to
T_{\gamma_{u}(t)}M$ is the parallel translation with respect to
$\nabla$ along $\gamma_{u}$ from $o=\gamma_{u}(0)$ to
$\gamma_{u}(t).$ Then $R_{u}(t)$ is a curve in the
$n(n+1)/2$-dimensional vector space ${\mathcal S}({\mathfrak m})$
of all self-adjoint operators of $({\mathfrak m},<,>)$ with
$R_{u}(0) = R_{u},$ being $R_{u}:= R(u,\cdot)u.$ ${\mathcal
S}({\mathfrak m})$ is a Euclidean vector space with inner product
defined by
\begin{equation}\label{inner}
<K,K'> = \sum_{i=1}^{n}<K(e_{i})K'(e_{i})>,
\end{equation}
for all $K,K'\in {\mathcal S}({\mathfrak m}),$ where
$\{e_{1},\dots,e_{n}\}$ is an arbitrary orthonormal basis of
$({\mathfrak m},<,>).$ Because, for naturally reductive spaces,
the {\em connection function} $\Lambda:{\mathfrak
m}\times{\mathfrak m}\to {\mathfrak m}$ associated to $\nabla$ is
given by $\Lambda_{x}(y) = 1/2[x,y]_{\mathfrak m} = S_{x}y,$ for
all $x,y\in {\mathfrak m},$ \cite[Theorem 2.10, p.197]{KN} it
directly follows
\[
\frac{\nabla (\exp tu)_{*o}v}{dt}_{\mid t} = (\exp tu)_{*o}S_{u}v,
\]
for each $v\in {\mathfrak m},$ and hence the parallel translation
$\tau_{u}$ is given by
\begin{equation}\label{trans}
\tau_{u}^{t} = (\exp tu)_{*o}e^{-tS_{u}},
\end{equation}
where $e$ denotes the exponential map of the Lie group of the
automorphisms $\Aut({\mathfrak m})$ of ${\mathfrak m}.$ Note that
$S_{u}:{\mathfrak m}\to {\mathfrak m}$ is a skew-symmetric
endomorphism of $({\mathfrak m},<,>)$ and so, $e^{S_{u}}$ is a
linear isometry of $({\mathfrak m},<,>).$
\begin{lemma}\label{lR} We have:
\begin{equation}\label{Ru}
R_{u}(t) = Ad_{e^{tS_{u}}}R_{u}.
\end{equation}
\end{lemma}
\noindent{\sf Proof.} For $x,y\in {\mathfrak m},$ using
(\ref{trans}) one obtains
$$
\begin{array}{lcl}
<R_{u}(t)x,y> & = &
g_{\gamma_{u}(t)}(R_{\gamma_{u}(t)}\tau^{t}_{u}x,\tau^{t}_{u}y) =
g_{\gamma_{u}(t)}(R_{\gamma_{u}(t)}(\exp
tu)_{*o}e^{-tS_{u}}x,(\exp tu)_{*o}e^{-tS_{u}}y)\\[0.5pc]
 & = & g_{\gamma_{u}(t)}((\exp tu)_{*o}R_{u}e^{-tS_{u}}x,(\exp
 tu)_{*o}e^{-tS_{u}}y)=
 <R_{u}e^{-tS_{u}}x,e^{-tS_{u}}y>\\[0.5pc]
 & = & <e^{tS_{u}}\comp R_{u}\comp e^{-tS_{u}}(x),y> =
 <(Ad_{e^{tS_{u}}}R_{u})x,y>.
 \end{array}
$$
Then, we have proved (\ref{Ru}).\hfill $\Box$

Hence, the curve $R_{u}(t)$ in ${\mathcal S}({\mathfrak m})$ can
be expressed as the power series expression
\begin{equation}\label{onepa}
R_{u}(t) = e^{t{\rm ad} S_{u}}(R_{u}) =
\sum_{k=0}^{\infty}\frac{t^{k}}{k!}S^{k}_{u}\cdot R_{u},
\end{equation}
where $S_{u}$ acts as a derivation on the space of the
endomorphisms $\End({\mathfrak m})$ of ${\mathfrak m}.$ It proves
the following.
\begin{lemma}\label{lrankos} We have:
\[
R_{u}^{i)}(0) = S_{u}^{i}\cdot R_{u},\;\;\;\mbox{for all}\;\;i\in
{\mathbb N}.
\]
\end{lemma}
 A curve $\alpha: I\to M$ in an arbitrary manifold $M$ is said
that has {\em constant osculating rank} $r$ if for all $t\in I,$
its higher order derivatives $\alpha'(t),\dots ,\alpha^{r)}(t)$
are linearly independent and $\alpha'(t),\dots,\alpha^{r+1)}(t)$
are linearly dependent in $T_{\alpha(t)}M.$ Orbits of
one-parameter subgroups of a Lie group acting on $M$ are examples
of curves with constant osculating rank. Concretely, we have
\begin{lemma}\label{losc} Let $G$ be a Lie group acting on $M$
from the left and let $\alpha$ be the curve in $M$ given by
$\alpha(t) = (\exp tX)p,$ for some $X\in {\mathfrak g}$ and $p\in
M.$ Then,
\[
\alpha^{k)}(t) = (\exp tX)_{*p}\alpha^{k)}(0),
\]
for all $t\in \RR$ and $k\in {\mathbb N},$ and so it has constant
osculating rank.
\end{lemma}
\noindent{\sf Proof.} From (\ref{adjun1}), taking into account
that $\alpha$ is the integral curve through $p$ of the fundamental
vector field $X^{*},$ we obtain $\alpha'(t) = (\exp
tX)_{*p}\alpha'(0).$ Hence, it follows that $\alpha'(t+s) = (\exp
tX)_{*\alpha(s)}\alpha'(s)$ and then,
$$
\begin{array}{lcl}
\alpha''(t) & = & \frac{d}{dt}_{\mid t}(\exp tX)_{*p}\alpha'(0) =
\frac{d}{dt}_{\mid t}\frac{d}{ds}_{\mid s=0}(\exp(t+s)X)p =
\frac{d}{ds}_{\mid s=0}\frac{d}{dt}_{\mid t}\alpha(t+s)\\[0.5pc]
&=  & \frac{d}{ds}_{\mid s=0}\alpha'(t+s)= \frac{d}{ds}_{\mid
s=0}(\exp tX)_{*\alpha(s)}\alpha'(s) = (\exp tX)_{*p}\alpha''(0).
\end{array}
$$
For the general case, we just use induction. Suppose that
$\alpha^{i)}(t) = (\exp tX)_{*p}\alpha^{i)}(0),$ for any $i\leq
k-1$ and $k\geq 3.$ We obtain
$$
\begin{array}{lcl}
\alpha^{k)}(t) & = & \frac{d}{dt}_{\mid t}(\exp
tX)_{*p}\alpha^{k-1)}(0) = \frac{d}{dt}_{\mid t}\frac{d}{ds}_{\mid
s=0}(\exp tX)_{*\alpha(s)}(\exp sX)_{*p}\alpha^{k-2)}(0)\\[0.5pc]
& = & \frac{d}{ds}_{\mid s=0}\frac{d}{dt}_{\mid
t}(\exp(t+s)X)_{*p}\alpha^{k-2)}(0) = \frac{d}{ds}_{\mid
s=0}\frac{d}{dt}_{\mid t}\alpha^{k-2)}(t+s) = \frac{d}{ds}_{\mid
s=0}\alpha^{k-1)}(t+s)\\[0.5pc]
 &= & \frac{d}{ds}_{\mid s=0}(\exp
tX)_{*\alpha(s)}\alpha^{k-1)}(s) = (\exp tX)_{*p}\alpha^{k)}(0).
\end{array}
$$
\hfill $\Box$

Because $\{Ad_{e^{tS_{u}}}\mid t\in \RR\}$ is a one-parameter
subgroup of the isometry group of $({\mathcal S}({\mathfrak
m}),<,>),$ it follows from Lemmas \ref{lR} and \ref{losc} that it
has constant osculating rank, the {\em Jacobi osculating rank of
$\gamma_{u}.$} For each unit vector $u\in {\mathfrak m},$ denote
by ${\mathcal R}_{u}({\mathfrak m})$ the smallest subspace of
${\mathcal S}({\mathfrak m})$ such that $R_{u}$ and $S_{u}\cdot
{\mathcal R}_{u}({\mathfrak m})\subset {\mathcal R}_{u}({\mathfrak
m}).$ From Lemma \ref{lrankos}, ${\mathcal R}_{u}({\mathfrak m})$
is generated by $R_{u},$ $S_{u}\cdot R_{u},\dots ,S^{r}_{u}\cdot
R_{u}.$ Hence, $r = \dim {\mathcal R}_{u}({\mathfrak m}) -1$ or $r
= \dim{\mathcal R}_{u}({\mathfrak m})$ and then, taking into
account that $K(u) = 0,$ for each $K\in {\mathcal
R}_{u}({\mathfrak m}),$ we have
\[
{\rm rank}_{\rm osc}(u)\leq \dim {\mathcal R}_{u}({\mathfrak
m})\leq \frac{n(n-1)}{2}.
\]
\begin{remark}\label{go}{\rm A homogeneous Riemannian manifold $(M = G/K,g)$
is said to be a {\em g.o. space} if each geodesic starting at the
origin is an orbit of an one-parameter subgroup $(\exp tZ)$ of
$Z\in {\mathfrak g}.$ Naturally reductive spaces are g.o. spaces
but there is a large number of examples of g.o. spaces which are
not naturally reductive. In similar way as before, one obtains
that the parallel translation $\tau_{u}$ along the geodesic
$\gamma_{u}(t) = (\exp tZ)o$ on a g.o. space is given by
\[
\tau^{t}_{u} = (\exp tZ)_{*o}e^{-t\Lambda_{u}}
\]
and then the formulas (\ref{Ru}), (\ref{onepa}) and Lemma
\ref{lrankos} hold. Hence, the notion of Jacobi osculating rank
can be directly extended to g.o. spaces (see \cite{AM-N} for more
details and references).}
\end{remark}

When the Jacobi osculating rank does not depend on the choice of
the geodesic, the naturally reductive homogeneous manifold $(M,g)$
is said to have {\em constant Jacobi osculating rank.}

\begin{lemma}\label{psim} Any naturally reductive homogeneous manifold
$(M,g)$ of constant Jacobi osculating rank zero is locally
symmetric.
\end{lemma}
\noindent{\sf Proof.} For each unit vector $u$ in ${\mathfrak m},$
let $\lambda_{1},\dots ,\lambda_{n}$ be the eigenvalues and
$\{e_{1},\dots ,e_{n}\}$ the corresponding eigenvectors of the
operator $R_{u}.$ Because ${\rm rank}_{\rm osc}(u) = 0,$ it
follows that $R_{u}(t) = R_{u},$ for all $t\in \RR,$ and then,
$\tau_{u}\comp R_{u} = R_{\gamma_{u}}\comp \tau_{u}.$ Hence,
$R_{\gamma_{u}}E_{i} = \lambda_{i}E_{i},$ $i=1,\dots ,n,$ where
$E_{i}$ is the vector field along $\gamma_{u}$ given by $E_{i}(t)
= \tau^{t}_{u}(e_{i}).$ Then $\{E_{i},\dots,E_{n}\}$ becomes into
a parallel frame field of eigenvectors of the Jacobi operator
$R_{\gamma_{u}}$ with constant eigenvalues $\lambda_{i}$ and so,
$(M,g)$ must be locally symmetric \cite{BeVa}.\hfill $\Box$

Because symmetric spaces have constant Jacobi osculating rank
zero, then we have
\begin{theorem}\label{sim} A simply connected, naturally reductive
homogeneous manifold has constant Jacobi osculating rank zero if
and only if it is a symmetric space.
\end{theorem}
Next, we also consider the curve $\tilde{R}_{u}(t)$ in ${\mathcal
S}({\mathfrak m})$ obtained by the $\nabla$-parallel translation
of $\tilde{R}_{\gamma_{u}}$ along $\gamma_{u}.$ First, using
(\ref{RR}), we get
\begin{equation}\label{RS}
\tilde{R}_{\gamma_{u}} = R_{\gamma_{u}} + S^{2}_{\gamma_{u}}.
\end{equation}
\begin{lemma}\label{pcharac} The curves $R_{u}(t)$ and $\tilde{R}_{u}(t)$ have the same osculating
rank. Moreover, we have:
\begin{enumerate}
\item[{\rm (i)}] $\tilde{R}_{u}(t) = R_{u}(t) + S^{2}_{u}.$
\item[{\rm (ii)}] $R_{u}^{i)} (0) = \tilde{R}_{u}^{i)}(0) = S^{i}_{u}\cdot
\tilde{R}_{u},$ for all $i\in {\mathbb N}.$
\end{enumerate}
\end{lemma}
\noindent{\sf Proof.} Because $\frac{\nabla}{dt} =
\frac{\tilde{\nabla}}{dt} + S_{\gamma_{u}},$ one obtains
$\frac{\nabla S_{\gamma_{u}}}{dt} =
\frac{\tilde{\nabla}S_{\gamma_{u}}}{dt}$ and, using that
$\tilde{\nabla}S =0,$ the tensor field $S_{\gamma_{u}}$ is
$\nabla$-parallel along $\gamma_{u},$ or equivalently,
$\tau^{t}_{u}\comp S_{u} = S_{\gamma_{u}}\comp \tau^{t}_{u}.$
Then, the curve $S_{u}(t) = \tau^{-t}_{u}\comp S_{\gamma_{u}}\comp
\tau^{t}_{u}$ is constant and (i) follows by using (\ref{RS}).
Now, using Lemma \ref{lrankos}, we also get (ii).\hfill $\Box$

 On naturally reductive homogeneous spaces $(M = G/K,g),$ any
$G$-invariant (unit) vector field is Killing and so, each one of
its integral curves is a geodesic.

\begin{lemma}\label{lnew} The Jacobi osculating rank of
each integral curve of a $G$-invariant vector field is zero.
\end{lemma}
\noindent{\sf Proof.} Put $u =U_{o}\in {\mathfrak m}.$ Then $u$ is
$Ad(K)$-invariant and it implies $[{\mathfrak k},u] = 0$ and from
(\ref{torsion}), $\tilde{R}_{u} =0.$ Because $\tilde{R}$ is
$G$-invariant, it follows that $\tilde{R}_{\gamma_{u}}\comp (\exp
tu)_{*o} = (\exp tu)_{*o}\comp\tilde{R}_{u}$ and then
$\tilde{R}_{\gamma_{u}}$ vanishes along $\gamma_{u}.$ From Lemma
\ref{pcharac}, ${\rm rank}_{\rm osc}(u) = 0.$ \hfill $\Box$

Hence, using Lemma \ref{psim}, we have

\begin{proposition}\label{nonconst} Any non-locally symmetric
naturally reductive space $(M = G/K,g)$ with a $G$-invariant
vector field has non-constant Jacobi osculating rank.
\end{proposition}

In \cite[Lemma 5.5]{GS} it is proved that any non-locally
symmetric naturally reductive space $(M,g)$ of dimension $n\leq 5$
admits a naturally reductive quotient representation $G/K$ and a
(non-parallel) $G$-invariant unit vector field. Then, we can
conclude
\begin{corollary} Any non-locally symmetric naturally reductive
space of dimension $n\leq 5$ has non-constant Jacobi osculating
rank.
\end{corollary}

\begin{remark}{\rm Simply connected Killing-transversally
symmetric spaces are introduced in \cite{GGV1} as simply connected
Riemannian manifolds equipped with a complete unit Killing vector
field $\xi$ such that all reflections with respect to its integral
curves are isometries. This family of naturally reductive spaces
contains to $M^{3}(\kappa,\tau),$ for all $\kappa$ and $\tau>0,$
and $\xi$ is $G$-invariant with respect to a naturally reductive
representation $M=G/K.$ The dual one-form of $\xi$ with respect to
the metric is a contact form if and only if it is irreducible
\cite[Theorem 5.1]{GGV1}. Then we can give a lot of examples of
irreducible naturally reductive spaces with non-constant Jacobi
osculating rank: Generalized Heisenberg groups equipped with
suitable left-invariant metrics, Berger's spheres,
$\varphi$-symmetric spaces, Sasakian space forms, etc.}
\end{remark}

\section{Jacobi osculating rank of geodesics on $M^{3}(\kappa,\tau)$}\indent

A naturally reductive decomposition for the Lie algebra of the
isometry group and an adapted canonical connection for
$M^{3}(\kappa,\tau)$ are obtained in \cite{TV} by using of
naturally reductive homogeneous structures. Next, we give a brief
summary from Theorems 6.4 and 6.5 in \cite{TV}. A non-vanishing
naturally reductive homogeneous structure on an arbitrary
three-dimensional (oriented) Riemannian manifold $(M^{3},g)$ can
be expressed as $S = \lambda dv,$ for some non-zero constant
$\lambda,$ where $dv$ is the volume form on $M^{3}$ and $S$ also
denotes the $3$-form $S(X,Y,Z) = g(S_{X}Y,Z),$ for all $X,Y,Z\in
{\mathfrak X}(M^{3}).$ Because $S$ and $-S$ are isomorphic
structures, we can take $\lambda>0.$ Then there exists an
orthonormal basis $\{e_{1}, e_{2},e_{3}\}$ of ${\mathfrak m}\cong
T_{o}M^{3}$ such that
\begin{equation}\label{T1}
S_{e_{1}}e_{2} = \lambda e_{3}, \;\;\; S_{e_{3}}e_{1} = \lambda
e_{2}, \;\;\; S_{e_{2}}e_{3} = \lambda e_{1}.
\end{equation}
For $M^{3}(\kappa,\tau),$ putting $e_{3} = \xi_{o},$ we get
$\lambda =\frac{\tau}{2}$ and $\tilde{R}$ is expressed as
\begin{equation}\label{R1}
\tilde{R}_{e_{1}e_{2}} = (\kappa - \tau^{2})A_{12},\;\;\;\;
\tilde{R}_{e_{1}e_{3}} = 0,\;\;\; \tilde{R}_{e_{2}e_{3}} = 0,
\end{equation}
where $A_{12}$ is the skew-symmetric endomorphism on ${\mathfrak
m}$ given by
\[
A_{12}e_{1} = e_{2},\;\;\;\; A_{12}e_{2} = -e_{1},\;\;\;\;
A_{12}e_{3} = 0.
\]
Note that if $\kappa=\tau^{2},$ then $\tilde{R}$ vanishes and
$M^{3}(\kappa,\tau)$ is an Einstein manifold and hence, of
constant curvature. The holonomy algebra ${\mathfrak k}$ of
$\tilde{\nabla}$ is generated by $A_{12}$ and the transvection
algebra ${\mathfrak tr}({\mathfrak m}) = {\mathfrak m}\oplus
{\mathfrak k}$ is generated by $\{e_{1},e_{2},e_{3},A_{12}\}$ with
Lie bracket given by
\begin{equation}\label{brac}
\left\{
\begin{array}{lcllcl}
[e_{1},e_{2}] =  \tau e_{3} + (\kappa-\tau^{2})A_{12}, & &
[e_{1},e_{3}]
 =  -\tau e_{2}, & & [e_{2},e_{3}]  =
\tau e_{1}\\[0.5pc]
[A_{12},e_{1}]  =  e_{2}, & & [A_{12},e_{2}]  = -e_{1}, & &
[A_{12},e_{3}]  = 0.
\end{array}
\right.
\end{equation}
Hence, putting
\begin{equation}\label{basen}
u_{1} = e_{1},\;\;\;\; u_{2} = e_{2},\;\;\;\; u_{3} = e_{3} +
\frac{\kappa-\tau^{2}}{\tau} A_{12},
\end{equation}
the subspace ${\mathfrak h}$ of ${\mathfrak tr}({\mathfrak m})$
generated by $u_{1},$ $u_{2},$ $u_{3}$ is a three-dimensional
unimodular Lie algebra and the linear isometry $f:
T_{o}M^{3}(\kappa,\tau)\cong{\mathfrak m} \to {\mathfrak h}$ given
by $f(e_{i}) = u_{i},$ $i=1,2,3,$ determines a (unique) isometry
from $M^{3}(\kappa,\tau)$ to the connected and simply connected
Lie group with Lie algebra ${\mathfrak h}$ equipped with the left
invariant metric such that $\{u_{1},u_{2},u_{3}\}$ is an
orthonormal basis. Then we can identify $M^{3}(\kappa,\tau)$ with
this unimodular Lie group and $u_{3}$ with the Hopf vector field
$\xi.$ From (\ref{basen}), we get
\begin{equation}\label{brac1}
\begin{array}{lcllcllcllcllcl}
[u_{1},u_{2}] & = & \tau u_{3}, & & [u_{2},u_{3}] & = &
\frac{\kappa}{\tau}u_{1}, & & [u_{3},u_{1}]&  = &
\frac{\kappa}{\tau}u_{2},\\[0.5pc]
[A_{12},u_{1}] & = & u_{2}, & & [A_{12},u_{2}] & = &-u_{1}, & &
[A_{12},u_{3}] & = & 0.
\end{array}
\end{equation}
Hence, it follows that ${\mathfrak tr}({\mathfrak m})$ is a
semi-direct sum ${\mathfrak h}\times_{\beta}{\mathfrak k}$ of
${\mathfrak h}$ and ${\mathfrak k}$ where $\beta$ is the
homomorphism of ${\mathfrak k}$ into $End({\mathfrak h})$ given by
$\beta(A_{12})u_{1} = u_{2},$ $\beta(A_{12})u_{2} = -u_{1},$
$\beta(A_{12})u_{3} = 0$ and we have
\[
\nabla_{X}\xi = \frac{\tau}{2}X\times \xi,
\]
where $\times$ denotes the vector product in ${\mathfrak h}.$ It
implies that $\xi$ is a unit Killing vector field and so, the
geodesic $\gamma_{e_{3}}(t)= (\exp te_{3})o$ coincides with its
integral curve through the origin. Moreover, the sectional
curvature $K(X,\xi)$ of the two-plane spanned by $X,$ $\xi$ is a
non-negative constant $c^{2},$ called the $\xi$-{\em sectional
curvature} \cite{GGV1}, given by $c^{2} = \frac{\tau^{2}}{4}$ and
$\{u_{1},u_{2},\xi\}$ is in fact a basis of eigenvectors for the
Ricci tensor $\rho$ with $\rho(u_{1},u_{1}) = \rho(u_{2},u_{2}) =
\kappa - \frac{\tau^{2}}{2}$ and $\rho(\xi,\xi)
=\frac{\tau^{2}}{2}.$ Furthermore, $\xi$ is the vertical field of
the fibration $M^{3}(\kappa,\tau) \to M^{3}/\xi = M^{2}(\kappa).$
Then $M^{3}(\kappa,\tau)$ is a principal $G^{1}$-bundle, where
$G^{1}$ denotes the one-parameter subgroup of global isometries
$\varphi_{t}$ generated by $\xi.$ $G^{1}$ is isomorphic to either
the circle group $S^{1}$ or to $\RR$ depending on whether the
integral curves of $\xi$ are closed or not. In particular, $G^{1}$
is a circle when $M^{3}(\kappa,\tau)$ is compact. If $c= 0,$ or
equivalently $\tau$ vanishes, this fibration becomes trivial and
we get the product spaces $M^{2}(\kappa) \times \RR.$ For $\kappa
>0,$ the metrics $g = g_{\kappa,\tau}$ are known as {\em Berger
metrics} and the corresponding fibration, as the {\em Hopf
fibration}. Here, as it has been already said, $\xi$ is called the
{\em Hopf vector field} of $M^{3}(\kappa,\tau),$ even for the
cases $\kappa \leq 0.$

We shall use the following result, which can be of interest by
itself.
\begin{lemma}\label{coclass} Each non-symmetric connected and simply connected
three-dimensional normal homogeneous manifold is isometric to
$S^{3}(\kappa,\tau),$ for some pair $(\kappa,\tau)$ such that
$\kappa >\tau^{2}.$
\end{lemma}
\noindent{\sf Proof.} The set of all inner products on ${\mathfrak
t}{\mathfrak r}({\mathfrak m})$ whose restriction to ${\mathfrak
m}$ is $<,>$ and such that ${\mathfrak k}$ is orthogonal to
${\mathfrak m}$ form an one-parameter family $\{B_{r}\mid r>0\},$
where $B_{r}(A_{12},A_{12}) = r.$ Then, from (\ref{tt2}) and using
(\ref{brac}), $B_{r}$ is bi-invariant if and only if $r=
\frac{1}{\kappa-\tau^{2}}.$ So, the existence of bi-invariant
metrics on $M^{3}(\kappa, \tau)$ is determined by the condition
$\kappa -\tau^{2}>0.$ \hfill $\Box$

Each $u\in {\mathfrak m}$ can be written as $u(\theta,\phi) =
\sin\theta\cos\phi e_{1} + \sin\theta \sin\phi e_{2} + \cos\theta
e_{3},$ where $\theta\in [0,\pi]$ and $\phi\in [0,2\pi].$ Then, we
get
\begin{equation}\label{A12}
u(\theta,\phi) = e^{\phi A_{12}}u(\theta),
\end{equation}
where $u(\theta)$ denotes the unit vector in the plane
$\RR\{e_{1},e_{3}\}$ given by $u(\theta) =u(\theta,0).$
\begin{lemma}\label{lad}We have:
\begin{enumerate}
\item[{\rm (a)}] $(\exp \phi A_{12})\gamma_{u(\theta)} =
\gamma_{u(\theta,\phi)}.$
\item[{\rm (b)}] $S_{u(\theta,\phi)} = Ad_{e^{\phi A_{12}}}S_{u(\theta)},$ $\tilde{R}_{u(\theta,\phi)} = Ad_{e^{\phi
A_{12}}}\tilde{R}_{u(\theta)},$ $R_{u(\theta,\phi)} = Ad_{e^{\phi
A_{12}}}R_{u(\theta)}.$
\item[{\rm (c)}] $R_{u(\theta,\phi)}(t) = Ad_{e^{\phi
A_{12}}}R_{u(\theta)}(t).$
\item[{\rm (d)}] ${\rm rank}_{\rm osc}u(\theta,\phi) = {\rm rank}_{\rm
osc}u(\theta),$ for all $\phi\in [0,2\pi].$
\end{enumerate}
\end{lemma}
\noindent {\sf Proof.} From (\ref{brac}), we get $Ad_{\exp\phi
A_{12}} = e^{\phi ad A_{12}}= e^{\phi A_{12}}.$ Then (\ref{A12})
can be written as
\[
u(\theta,\phi) = Ad_{\exp \phi A_{12}}u(\theta).
\]
Hence, (a) follows taking into account that the action of the
linear isotropy group of the isotropy subgroup $K$ of
$M^{3}(\kappa,\tau)$ at the origin corresponds under projection
with $Ad(K)$ on ${\mathfrak m}.$ Note that $K$ is connected
because $M^{3}(\kappa,\tau)$ is simply connected and then it
coincides with the one-parameter subgroup $\exp \phi A_{12}.$

In (b), we shall only prove the first equality, the other two are
obtained in similar way. Because the tensor field $S$ is invariant
for the isometry group of $M^{3}(\kappa,\tau),$ $S$ as tensor on
${\mathfrak m}$ is $Ad(K)$-invariant and then $S_{u(\theta,\phi)}
\comp Ad_{\exp \phi A_{12}} = Ad_{\exp\phi A_{12}}\comp
S_{u(\theta)}.$ Hence, we obtain
\[
S_{u(\theta,\phi)} = Ad_{Ad_{\exp\phi A_{12}}}S_{u(\theta)} =
Ad_{e^{\phi A_{12}}}S_{u(\theta)}.
\]
For (c), we use Lemma \ref{lrankos} and (b) to obtain
\begin{equation}\label{Rad}
R^{i)}_{u(\theta,\phi)}(0) = Ad_{e^{\phi
A_{12}}}R^{i)}_{u(\theta)}(0).
\end{equation}
Finally, because $Ad_{e^{\phi A_{12}}}$ belongs to $\Aut({\mathcal
S}({\mathfrak m})),$ we have (d).\hfill $\Box$

Therefore, we may restrict our study of the geometry of geodesics
on $M^{3}(\kappa,\tau)$ to geodesics emanating from the origin
with initial directions $u(\theta),$ $\theta\in [0,\pi].$

\begin{theorem}\label{tosc} The Jacobi osculating rank of every geodesic $\gamma_{u}$ on $M^{3}
(\kappa,\tau)$ is two except for the Hopf fibers, where it is
zero. Moreover, we have:
\begin{enumerate}
\item[{\rm (i)}] $R_{u}(t) = R_{u} + \frac{1}{\tau}\Big (
\sin\tau t R'_{u}(0) + \frac{1}{\tau}(1-\cos\tau t)R''_{u}(0)\Big
).$
\item[{\rm (ii)}] For each $u = u(\theta,\phi)$ with $\theta\in ]0,\pi[,$ $R_{u}(t)$ is a circle in ${\mathcal S}({\mathfrak m})$ of radius
$\frac{\sqrt{2}}{2}|\tau^{2}-\kappa|\sin^{2}\theta.$
\end{enumerate}
\end{theorem}

\noindent {\sf Proof}. With respect to the basis
$\{e_{1},e_{2},e_{3}\}$ on ${\mathfrak m},$ using (\ref{T1}) and
(\ref{R1}), we get
$$
S_{u(\theta)} = \frac{\tau}{2} \left (
\begin{array}{ccc}
0 & -\cos\theta & 0\\[0.4pc]
\cos\theta & 0 & -\sin\theta\\[0.4pc]
0 & \sin\theta & 0
\end{array}
\right ) ; \;\;\;\;\; \tilde{R}_{u(\theta)} =
2\tau^{2}\mu(\theta)\left (
\begin{array}{ccc}
0 & 0 & 0\\[0.4pc]
0 & -1 & 0\\[0.4pc]
0 & 0 & 0
\end{array}
\right ),
$$
where $\mu(\theta)=
\frac{\tau^{2}-\kappa}{2\tau^{2}}\sin^{2}\theta.$ Note that
$\mu(\theta) = 0$ if and only if $\theta = 0$ or $\theta = \pi,$
i.e. on the Hopf direction. By a direct computation, we have
$$
S_{u(\theta)}\cdot \tilde{R}_{u(\theta)} =
\tau^{3}\mu(\theta)\left (
\begin{array}{ccc}
0 & \cos\theta & 0\\[0.4pc]
\cos\theta & 0 & -\sin\theta\\[0.4pc]
0 & -\sin\theta & 0
\end{array}
\right )
$$
and
$$
S_{u(\theta)}^{2}\cdot \tilde{R}_{u(\theta)} =
\tau^{4}\mu(\theta)\left (
\begin{array}{ccc}
-\cos^{2}\theta & 0 & \sin\theta\cos\theta\\[0.4pc]
0 & 1 & 0\\[0.4pc]
\sin\theta\cos\theta & 0 & -\sin^{2}\theta
\end{array}
\right ) .
$$
Moreover, one obtains $S^{3}_{u(\theta)}\cdot
\tilde{R}_{u(\theta)}  =-\tau^{2}S_{u(\theta)}\cdot
\tilde{R}_{u(\theta)}$ and $S^{4}_{u(\theta)}\cdot
\tilde{R}_{u(\theta)}  = -\tau^{2}S_{u(\theta)}^{2}\cdot
\tilde{R}_{u(\theta)}.$ Therefore, taking into account that
$S^{k+1}\cdot\tilde{R}_{u(\theta)} = S_{u(\theta)}\comp(S^{k}\cdot
\tilde{R}_{u(\theta)}) - (S^{k}\cdot\tilde{R}_{u(\theta)})\comp
S_{u(\theta)},$ it follows by the induction
\begin{equation}\label{induc}
\begin{array}{lcl}
S^{2k-1}_{u(\theta)}\cdot \tilde{R}_{u(\theta)} & = &
(-1)^{k-1}\tau^{2(k-1)}S_{u(\theta)}\cdot \tilde{R}_{u(\theta)},\\[0.5pc]
S^{2k}_{u(\theta)}\cdot \tilde{R}_{u(\theta)} & = &
(-1)^{k-1}\tau^{2(k-1)} S_{u(\theta)}^{2}\cdot
\tilde{R}_{u(\theta)}.
\end{array}
\end{equation}
Then, from Lemma \ref{pcharac} and Lemma \ref{lad} (d), we obtain
that ${\rm rank}_{\rm osc}u(\theta,\phi) = 2$ if $\theta\in
]0,\pi[$ and ${\rm rank}_{\rm osc}(\pm e_{3}) = 0$ and it proves
the first part of the Theorem. Moreover, (i) follows from
(\ref{induc}), using Lemma \ref{pcharac}, Lemma \ref{lad} (b) and
(\ref{Rad}).

For $\theta\in ]0,\pi[,$ the elements
$\{v_{1}(\theta),v_{2}(\theta),v_{3}(\theta)\}$ of ${\mathcal
S}({\mathfrak m})$ given by the matrices
$$
\begin{array}{lcl}
v_{1}(\theta) & = & \frac{\sqrt{2}}{2} \left (
\begin{array}{ccc}
0 & \cos\theta & 0\\[0.4pc]
\cos\theta & 0 & -\sin\theta\\[0.4pc]
0 & -\sin\theta & 0
\end{array}
\right ),\;\;\; v_{2}(\theta) = \frac{\sqrt{2}}{2} \left (
\begin{array}{ccc}
-\cos^{2}\theta & 0 & \sin\theta\cos\theta\\[0.4pc]
0 & 1 & 0\\[0.4pc]
\sin\theta\cos\theta & 0 & -\sin^{2}\theta
\end{array}
\right ) ,\\[2.5pc]
v_{3}(\theta) & = & \frac{\sqrt{2}}{2} \left (
\begin{array}{ccc}
-\cos^{2}\theta & 0 & \sin\theta\cos\theta\\[0.4pc]
0 & -1 & 0\\[0.4pc]
\sin\theta\cos\theta & 0 & -\sin^{2}\theta
\end{array}
\right )
\end{array}
$$
constitute an orthonormal basis for $R_{u(\theta)}({\mathfrak m})$
in $({\mathcal S}({\mathfrak m}),<,>)$ and we get
\begin{equation}\label{RRR}
\begin{array}{lcllcllcl}
\tilde{R}_{u(\theta)}\hspace{-0.20cm} &\hspace{-0.20cm} =
\hspace{-0.20cm}&
\hspace{-0.20cm}\sqrt{2}\tau^{2}\mu(\theta)(v_{3}(\theta) -
v_{2}(\theta)), & & S^{2}_{u(\theta)}
\hspace{-0.20cm}&\hspace{-0.20cm} =\hspace{-0.20cm}
& \hspace{-0.20cm}\frac{\sqrt{2}}{4}\tau^{2}v_{3}(\theta),\\[0.5pc]
R'_{u(\theta)}(0) \hspace{-0.20cm}&\hspace{-0.20cm}
=\hspace{-0.20cm} &\hspace{-0.20cm}
\sqrt{2}\tau^{3}\mu(\theta)v_{1}(\theta), & &
R''_{u(\theta)}(0)\hspace{-0.20cm} & \hspace{-0.20cm} =
\hspace{-0.20cm}& \hspace{-0.20cm}\sqrt{2}\tau^{4}\mu(\theta)v_{2}(\theta),\\[0.5pc]
R_{u(\theta)}\hspace{-0.20cm} &\hspace{-0.20cm} =\hspace{-0.20cm}
&\hspace{-0.20cm} -\frac{\sqrt{2}}{4}\tau^{2}\Big
(4\mu(\theta)v_{2}(\theta) + (1 - 4\mu(\theta))v_{3}(\theta)\Big
).
\end{array}
\end{equation}
Hence, (i) implies that
\[
R_{u(\theta)}(t) = \sqrt{2}\tau^{2}\mu(\theta)\Big (\sin \tau
tv_{1}(\theta) - \cos \tau tv_{2}(\theta)\Big ) +
\frac{\sqrt{2}}{4}\tau^{2}(\mu(\theta) - 1)v_{3}(\theta).
\]
Then, putting $R_{u(\theta)}({\mathfrak m})\cong \RR^{3}[x,y,z],$
where $x,$ $y,$ $z$ are the cartesian coordinates with respect to
$\{v_{1}(\theta),v_{2}(\theta),v_{3}(\theta)\},$
$R_{u(\theta)}(t)$ is the circle in the plane
$z=\frac{\sqrt{2}}{4}\tau^{2}(\mu(\theta)-1)$ such that $x=
\sqrt{2}\tau^{2}\mu(\theta)\sin\tau t,$ $y =
-\sqrt{2}\tau^{2}\mu(\theta)\cos\tau t.$ Because the Euclidean
product $<,>$ of ${\mathcal S}({\mathfrak m}),$ defined in
(\ref{inner}), coincides with the trace form of ${\rm
End}({\mathfrak m}),$ it is, in particular, $Ad(e^{\phi
A_{12}})$-invariant and so, (ii) follows from Lemma \ref{lad} (c),
taking into account that $S^{2}_{u(\theta,\phi)} = Ad_{e^{\phi
A_{12}}}S^{2}_{u(\theta)}.$ \hfill $\Box$

\begin{remark}{\rm The circle $R_{u(\theta,\phi)}(t)$ is in the plane of
${\mathcal S}({\mathfrak m})$ determined by its centre, i.e. the
point $(4\mu(\theta)-1)S^{2}_{u(\theta,\phi)},$ and the subspace
generated by $R'_{u(\theta,\phi)}(0)$ and
$R''_{u(\theta,\phi)}(0);$ its period is $\frac{2\pi}{\tau}$ and
it contracts to $R_{\xi_{0}}$ when $\theta$ converges to $0$ or to
$\pi.$}
\end{remark}

\section{Isotropic conjugate points in $M^{3}(\kappa,\tau)$}\indent

Next, we determine all pairs of conjugate points of
$M^{3}(\kappa,\tau)$ and those which are isotropic or strictly
isotropic. As before, we only need to consider geodesics
$\gamma_{u(\theta)},$ with $\theta\in [0,\pi].$ Note that if
$\gamma_{u(\theta)}$ is an isotropic geodesic, then
$\gamma_{u(\theta,\phi)}$ is also isotropic. From (\ref{ini}), a
Jacobi field $V$ along $\gamma_{u(\theta)}$ in
$M^{3}(\kappa,\tau),$ with $V(0) = 0$ is isotropic if and only if
$V'(0)$ is collinear with $e_{2}.$ Hence, $\dim {\rm
Isot}(\gamma_{u(\theta,\phi)}) = 1$ if $\theta\in ]0,\pi[$ and
zero for the fibers of $M^{3}(\kappa,\tau),$ where ${\rm
Isot}(\gamma_{u(\theta,\phi)})$ denotes the space of all isotropic
Jacobi fields along $\gamma_{u(\theta,\phi)}$ vanishing at the
origin. Moreover, $V$ along $\gamma_{u(\theta,\phi)}$ is isotropic
if and only if $V'(0)$ is collinear with $e^{\phi A_{12}}e_{2} =
-\sin \phi e_{1} + \cos \phi e_{2}.$ Therefore, we also have
\begin{lemma}\label{stric} All pair of isotropic conjugate points in $M^{3}(\kappa,\tau)
$ are strictly isotropic and all isotropic geodesic is strictly
isotropic.
\end{lemma}
Now, we prove the main theorem of this section.
\begin{theorem}\label{pral1} A geodesic $\gamma(t)$ on $M^{3}(\kappa,\tau)$ starting at the origin with
slope angle $\theta$ admits conjugate points to the origin if and
only if $\lambda(\theta) > 0,$ where $\lambda$ is the function
given by $\lambda(\theta): = \kappa \sin^{2}\theta +
\tau^{2}\cos^{2}\theta.$ Moreover, we have:
\begin{enumerate}
\item[{\rm (i)}] The conjugate points along the Hopf fibers
are at $t = \frac{2\pi p}{\tau},$ $p\in {\mathbb N},$ their
multiplicity is $2$ and they are not isotropic.

\item[{\rm (ii)}] For $\theta\in ]0,\pi[,$ the
conjugate points along $\gamma$ to the origin are all
$\gamma(\frac{s}{\sqrt{\lambda(\theta)}}),$ where
\begin{enumerate}
\item[{\rm 1.}] $s= 2 p\pi,$ $p\in {\Bbb N},$

\noindent or

\item[{\rm 2.}] $s$ is a solution of the equation $\tan \frac{s}{2}
= \mu(\theta)s,$ where $\mu(\theta) =
\frac{\tau^{2}-\kappa}{2\tau^{2}}\sin^{2}\theta.$
\end{enumerate}
In the first case, they are strictly isotropic and the second one,
they are not isotropic. In both cases, their multiplicity is $1.$
\end{enumerate}
\end{theorem}
\noindent{\sf Proof}. First, we shall show that the {\em tangent
conjugate locus} ${\rm conj}(M^{3}(\kappa,\tau))$ of the origin is
given by
\[
{\rm conj}(M^{3}(\kappa,\tau)) =
\{\frac{s}{\sqrt{\lambda(\theta)}}u(\theta,\phi)\mid
\lambda(\theta)>0, \;s\in{\mathcal Z }^{+}(f_{\theta})\},
\]
where ${\mathcal Z}^{+}(f_{\theta})$ is the set of zeros $s\in
{\RR}^{+}$ of $f_{\theta}(s) = 1-\cos s -\mu(\theta)s\sin s.$ (See
{\it Fig. 1} and {\it Fig. 2}).

The Jacobi equation (\ref{Jam1}) for geodesics
$\gamma_{u(\theta)}$ may be expressed, from (\ref{T1}) and
(\ref{R1}), as
\begin{equation}\label{sis}
\left \{
\begin{array}{l}
{X^{1}}'' - \tau\cos\theta{X^{2}}'=0,\\
{X^{2}}'' + \tau(\cos\theta{X^{1}}' - \sin\theta{X^{3}}')
-2\tau^{2}\mu(\theta)X^{2}=0,\\
{X^{3}}'' + \tau\sin\theta{X^{2}}' = 0,
\end{array}
\right.
\end{equation}
where the solutions $X(t)$ in ${\mathfrak m}$ are given by $X(t) =
\sum_{i=1}^{3}X^{i}(t)e_{i}.$ If $\theta = 0$ or $\theta= \pi,$
then it reduces to
$$
\left \{
\begin{array}{l}
{X^{1}}'' - \tau{X^{2}}' =0,\\
{X^{2}}'' + \tau{X^{1}}' = 0,\\
{X^{3}}'' = 0
\end{array}
\right.
$$
and the Jacobi solutions $X$ such that $X(0)=0$ are given by
$$
X(t)  =  (A(1-\cos \tau t) - B\sin\tau t)e_{1} + (A\sin\tau t +
B(1-\cos \tau t)e_{2} + Cte_{3},
$$
where $A,$ $B,$ $C$ are constant. Hence, the conjugate points to
the origin along the geodesic $\gamma(t) = (\exp te_{3})o$ are
given by $\gamma(2\pi p/\tau),$ for $p\in {\Bbb Z}^{*},$ the
multiplicity of each one of them is $2$ and they are not isotropic
using (\ref{ini}). It proves (i).

Next, we suppose that $\theta\in ]0,\pi[.$ Differentiating the
second equality and substituting in it the first and the second
one, the system (\ref{sis}) can be reduced to
\begin{equation}\label{sis2}
\left \{
\begin{array}{l}
{X^{1}}'' = \tau \cos \theta Y,\\
Y'' + \lambda Y = 0,\\
{X^{3}}''= -\tau\sin\theta Y,
\end{array}
\right.
\end{equation}
where $Y = {X^{2}}'.$ If $\lambda=\lambda(\theta)>0,$ then $Y =
A\cos\sqrt{\lambda}t + B \sin\sqrt{\lambda}t,$ where $A$ and $B$
are constants. It is straightforward to check that the solutions
$X^{i}(t)$ of above system such that $X^{i}(0) = 0,$ $i=1,2,3,$
for $\theta\neq \frac{\pi}{2}$ are given by
$$
\begin{array}{lcl}
X^{1}(t) & = &\frac{\tau\cos \theta}{\lambda}\Big
(A(1-\cos\sqrt{\lambda}t) - B\sin\sqrt{\lambda}t + C\sqrt{\lambda}t\Big ),\\
X^{2}(t) &  = & \frac{1}{\sqrt{\lambda}}\Big (A\sin\sqrt{\lambda}t
+B(1-\cos\sqrt{\lambda}t)\Big ),\\
X^{3}(t) & = & -\frac{\tau\sin \theta}{\lambda}\Big
(A(1-\cos\sqrt{\lambda}t) +
B(\frac{\tau^{2}-\kappa}{\tau^{2}}\sqrt{\lambda}t - \sin
\sqrt{\lambda}t) - C\cot^{2}\theta \sqrt{\lambda}t\Big ).
\end{array}
$$
Then the arc length $t,$ $t\neq 0,$ at the conjugate points along
$\gamma$ are the zeros of the determinant
$$
\left |
\begin{array}{ccc}
1-\cos\sqrt{\lambda}t & -\sin\sqrt{\lambda}t & 1\\[0.5pc]
\sin\sqrt{\lambda}t & 1 - \cos\sqrt{\lambda}t & 0\\[0.5pc]
1- \cos \sqrt{\lambda}t &
\frac{\tau^{2}-\kappa}{\tau^{2}}\sqrt{\lambda}t -
\sin\sqrt{\lambda}t & -\cot^{2}\theta
\end{array}
\right | .
$$
For $\theta =\frac{\pi}{2},$ the corresponding solutions are given
by
$$
\begin{array}{lcl}
X^{1}(t) & = & C t\\
X^{2}(t) &  = & \frac{1}{\sqrt{\kappa}}\Big (A\sin\sqrt{\kappa}t
+B(1-\cos\sqrt{\kappa}t)\Big ),\\
X^{3}(t) & = & -\frac{\tau}{\kappa}\Big (A(1-\cos\sqrt{\kappa}t) +
B(\frac{\tau^{2}-\kappa}{\tau^{2}}\sqrt{\kappa}t -
\sin\sqrt{\kappa} t\Big ).
\end{array}
$$

\noindent Hence, making straight calculations and putting $s =
s(t) = \sqrt{\lambda}t,$ one obtains that for both cases the
problem to find the conjugate points to the origin, reduces to
find the zeros of the function $f_{\theta}(s) = 1-\cos s -
\mu(\theta)s\sin s,$ for $\theta\in [0,\pi].$ Hence, it follows
that $\sin s=0,$ or equivalently $s\in 2\pi{\Bbb Z},$ or, putting
$\mu = \mu(\theta),$
\[
\cos s =\frac{1-\mu^{2}s^{2}}{1 + \mu^{2}s^{2}},\;\;\;\;\; \sin s
= \frac{2\mu s}{1 + \mu^{2}s^{2}},
\]
which yields to the equation $\tan \frac{s}{2} = \mu s.$ Because
the rank of these matrices is two, the multiplicity
$n_{\gamma_{u(\theta)}}(s/\sqrt{\lambda(\theta)})$ of the
conjugate point $\gamma_{u(\theta)}(s/\sqrt{\lambda(\theta)})$ is
one. Moreover, the space of isotropic Jacobi fields along $\gamma$
are spanned by $V(t) = (\exp tu)_{*o}X(t),$ where
\[
X(t)  =  \tau\cos\theta(1-\cos \sqrt{\lambda}t)e_{1} +
\sqrt{\lambda}\sin\sqrt{\lambda}t e_{2} - \tau\sin\theta(1-\cos
\sqrt{\lambda}t)e_{3}.
\]
Therefore, the isotropic conjugate points to the origin are all
$\gamma(t)$ with $t\sqrt{\lambda(\theta)}\in 2\pi {\Bbb Z}^{*}.$

For the rest of the cases, that is for $\lambda(\theta)\leq 0,$
one can check by straightforward computations that the geodesic
$\gamma_{u(\theta)}$ does not admit conjugate points to the
origin, which proves the first part of the theorem. We have also
proved that ${\mathcal Z }^{+}(f_{\theta})= \{2p\pi\mid
p\in{\mathbb N}\},$ if $\theta  = 0$ or $\theta = \pi$ and
${\mathcal Z }^{+}(f_{\theta})= \{2p\pi\mid p\in {\mathbb
N}\}\bigcup \{s\in {\RR}\mid \tan \frac{s}{2} = \mu(\theta)s\}$ if
$\theta\in ]0,\pi[.$ Hence, using Lemma \ref{stric}, we get (i)
and (ii). \hfill $\Box$

Next, we give some applications of Theorem \ref{pral1}. We start
determining the isotropic geodesics of $M^{3}(\kappa,\tau).$

\subsection{Isotropic geodesics} As well-known, every maximal geodesic of a homogeneous Riemannian
manifold is either one-to-one or simply closed.
\begin{proposition}\label{pl} We have:
\begin{enumerate}
\item[{\rm (i)}] All closed geodesic on $M^{3}(\kappa,\tau)$ starting at the origin with slope angle $\theta
\in ]0,\pi[$ admits isotropic conjugate points and its length is a
integer multiple of $\frac{2\pi}{\sqrt{\lambda(\theta)}}.$
\item[{\rm (ii)}] Any isotropic geodesic is one-to-one.
\end{enumerate}
\end{proposition}
\noindent{\sf Proof.} Let $\gamma_{u}$ be a closed geodesic
starting at the origin with length $l$ and slope angle $\theta\in
]0,\pi[.$ Then the vector field $V = A_{12}^{*}\comp \gamma_{u}$
along $\gamma_{u}$ satisfies $V(0) = V(l)=0.$ Moreover, it is an
non zero isotropic Jacobi field. In fact, we have (see, for
example \cite[p. 577]{Z})
\[
V'(0) = [A_{12},u] = A_{12}(u) \neq 0.
\]
Hence, $\gamma_{u}(0)$ and $\gamma_{u}(l)$
$(\gamma_{u}(0)=\gamma_{u}(l)$ is the origin of
$M^{3}(\kappa,\tau))$ are isotropic conjugate points and, using
Theorem \ref{pral1}, there exists $p\in {\mathbb N}$ such that $l=
\frac{2p\pi}{\sqrt{\lambda(\theta)}}.$ It proves (i). Because the
existence of isotropic conjugate points to the origin along a
geodesic implies also, using Theorem \ref{pral1}, the existence of
non-isotropic conjugate points, (i) allows to show (ii). \hfill
$\Box$

In \cite[Corollary 3.8 and Proposition 3.10]{Engel}, S. Engel has
proved that a geodesic on $M^{3}(\kappa,\tau)$ with slope angle
$\theta\in ]0,\pi[,$ such that $\lambda(\theta)>0$ intersects the
Hopf fiber through the origin periodically at the points
$\gamma(2p\pi/\sqrt{\lambda(\theta)}),$ $p\in {\mathbb N}.$ Then,
using Theorem \ref{pral1}, we can show the following.
\begin{proposition} Any geodesic on $M^{3}(\kappa,\tau)$ starting at the origin with slope angle $\theta\in
]0,\pi[$ intersects the Hopf fiber through the origin exactly at
its isotropic conjugate points.
\end{proposition}

\noindent From here and using Theorem \ref{pral1}, we have
\begin{corollary} Let $\gamma$ be a geodesic in
$M^{3}(\kappa,\tau)$ with slope angle $\theta\in ]0,\pi[.$ Then
the following statements are equivalent:
\begin{enumerate}
\item[{\rm (i)}] $\gamma$ is isotropic;
\item[{\rm (ii)}] $\gamma$ does not admit any pair of conjugate
points;
\item[{\rm (iii)}] $\gamma$ intersects each Hopf fiber only at a
unique point.
\end{enumerate}
\end{corollary}

\begin{remark}{\rm The Hopf fibers of $H_{3}(\tau)$ and
$\widetilde{SL_2}(\kappa,\tau)$ are one-to-one geodesics and on
$S^{3}(\kappa,\tau)$ they are simple closed geodesics with the
same length (see \cite{GGV1}). For the Berger metric obtained upon
scaling the fibers of the length $2\pi$ of the Hopf fibration
$S^{3}(1)\to S^{2}(4) = S^{3}/S^{1},$ where $S^{3}(\kappa)$ or
$S^{2}(\kappa)$ denote the correspoding spheres of radius
$\frac{1}{\sqrt{\kappa}},$ the length of its fibers is $l= |2\pi
c|,$ $c^{2}$ being the $\xi$-sectional curvature. For
$S^{3}(\kappa,\tau) \to S^{2}(\kappa),$ taking into account that
metric of $S^{2}(\kappa)$ is given by a homothetic change from the
metric of $S^{2}(4)$ with coefficient $\frac{4}{\kappa},$ it
follows that $l$ is given by $l = \frac{4\pi\tau}{\kappa}.$}
\end{remark}

Using the above results and taking into account that a geodesic
without pairs of conjugate points is considered by definition
strictly isotropic, we can conclude with the following
corollaries.
\begin{corollary}On $H_{3}(\tau),$ we have:
\begin{enumerate}
\item[{\rm (i)}] All geodesic is one-to-one.
\item[{\rm (ii)}] A geodesic is isotropic if and only if its slope
angle is $\theta= \pi/2.$
\end{enumerate}
\end{corollary}

\begin{corollary}On $\widetilde{SL_{2}}(\kappa,\tau),$ we have:
\begin{enumerate}
\item[{\rm (i)}] All geodesic is one-to-one.
\item[{\rm (ii)}] A geodesic is isotropic if and only if its slope
angle $\theta$ belongs to $[\varepsilon,  \pi- \varepsilon],$
where $\varepsilon = \arctan \frac{\tau}{\sqrt{-\kappa}}.$
\end{enumerate}
\end{corollary}

According with \cite[Corollary 3.12]{Engel}, a geodesic on
$S^{3}(\kappa,\tau)$ with slope angle $\theta\in ]0,\pi[$ is
closed if and only if
$\frac{\tau^{2}-\kappa}{\sqrt{\lambda(\theta)}}\cos\theta\in
{\mathbb Q}.$ Then, using Proposition \ref{pl}, it follows
\begin{corollary} On $S^{3}(\kappa,\tau),$ we have:
\begin{enumerate}
\item[{\rm (i)}] A geodesic with slope angle $\theta\in ]0,\pi[$
is closed if and only if its length $l$ satisfies
\[
\frac{l}{2\pi}(\tau^{2}-\kappa)\cos\theta\in {\mathbb Q}.
\]
\item[{\rm (ii)}] Any geodesic is not isotropic and it admits conjugate points.
\end{enumerate}
\end{corollary}

\subsection{Conjugate radius} Denote by $\rho_{\rm conj}(\theta)$ the {\em conjugate radius} of
a geodesic $\gamma$ on $M^{3}(\kappa,\tau)$ with slope angle
$\theta$ emanating from the origin. If $\rho_{\rm conj}(\theta)$
is finite, it is the distance along $\gamma$ from the origin to
its first conjugate point. Put $I_{M^{3}(\kappa,\tau)} =
\{\theta\in [0, \pi]\mid \lambda(\theta)>0\}.$ Then,
$I_{H_{3}(\tau)} = [0, \frac{\pi}{2}[\cup ] \frac{\pi}{2}, \pi],$
$I_{S^{3}(\kappa,\tau)} = [0,\pi]$ and
$I_{\widetilde{SL_{2}}(\kappa,\tau)} = [0,\varepsilon[\cup ]\pi
-\varepsilon, \pi],$ where $\varepsilon = \arctan
\frac{\tau}{\sqrt{-\kappa}}.$ Hence, one obtains $\mu(\theta)<
1/2$ for all $\theta\in I_{M^{3}(\kappa,\tau)}$ and moreover,
$0\leq \mu(\theta)<1/2,$ if $\kappa<\tau^{2},$ and
$\mu(\theta)\leq 0,$ if $\kappa> \tau^{2}.$ See {\it Fig. 1} and
{\it Fig. 2} for the behaviour of $f_{\theta}(s)$ in both cases.

\includegraphics{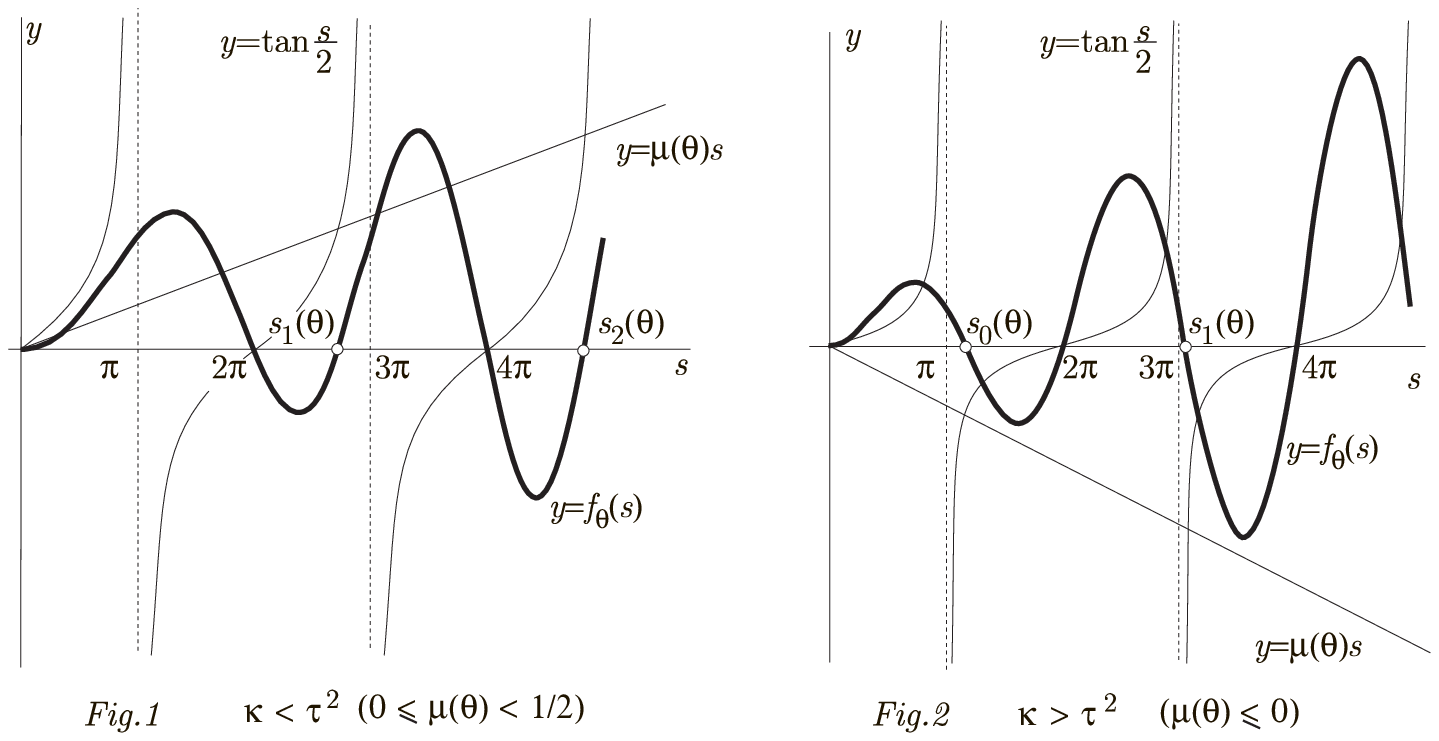}
\begin{proposition}\label{first} We have:
\begin{enumerate}
\item[{\rm (i)}] The conjugate radius $\rho_{\rm conj}(\theta)$ satisfies
$$
\rho_{\rm conj}(\theta) = \left\{
\begin{array}{l}
+\infty \;\;\; {\mbox if}\; \theta\notin
I_{M^{3}(\kappa,\tau)},\\[0.6pc]
 \frac{2\pi}{\sqrt{\lambda(\theta)}}\;\;\; {\mbox if}\;
\kappa<
\tau^{2},\; \theta\in I_{M^{3}(\kappa,\tau)}, \\[0.6pc]
\frac{s_{0}(\theta)}{\sqrt{\lambda(\theta)}},\;\;\; {\mbox if}\;
\kappa> \tau^{2},
\end{array}
\right.
$$
where $s_{0}(\theta)\in ]\pi,2\pi[$ and $\tan
\frac{s_{0}(\theta)}{2} = \mu(\theta)s_{0}(\theta).$

\item[{\rm (ii)}] The first conjugate point of $\gamma$ is
isotropic if and only if $\kappa <\tau^{2}$ and $\theta\in
]0,\pi[.$
\end{enumerate}
\end{proposition}
\noindent{\sf Proof.} Because $(\tan \frac{s}{2})'(0) =
\frac{1}{2},$ the equation $\tan \frac{s}{2} = \mu(\theta)s$ has
solutions in $]0, 2\pi[$ if and only if $\mu(\theta)\in]-\infty,
0[\cup ]1/2, +\infty [.$ Then, from Theorem \ref{pral1},
$\rho_{\rm conj}(\theta) = {2\pi}/\sqrt{\lambda(\theta)}$ if
$0\leq \mu(\theta) < 1/2$ and $\rho_{\rm conj}(\theta) =
s_{0}(\theta)/\sqrt{\lambda(\theta)}$ if $\mu(\theta)\in ]-\infty
, 0[.$ It gives the result. \hfill $\Box$

\begin{remark}{\rm The first conjugate points to the origin for the Berger's spheres
$M^{3}_{\alpha}$ given by $M^{3}_{\alpha} = S^{3}(4,2\sin\alpha),$
$\alpha\in]0,\pi/2],$ has been determined by I. Chavel \cite{Ch2}
(see also \cite[Lemma 2.3]{Sakai}, \cite[Proposition 3.1]{Ra}).
This result is a particular case of the above proposition.}
\end{remark}

For $\kappa > \tau^{2},$ one obtains that $\theta = \pi/2$
minimises $\mu(\theta)$ and so it is also a minimum for the map
$s_{0}.$ Hence, the following result is immediate.

\begin{corollary} The conjugate radius $\rho_{\rm conj}(M^{3}(\kappa,\tau))$ of
$M^{3}(\kappa,\tau)$ satisfies
$$
\rho_{conj}(M^{3}(\kappa,\tau)) = \left\{
\begin{array}{l}
\frac{2\pi}{\tau},\;\;\; {\mbox if}\; \kappa<
\tau^{2}\\[0.6pc]
\frac{s_{0}(\frac{\pi}{2})}{\sqrt{\kappa}},\;\;\; {\mbox if}\;
\kappa> \tau^{2}.
\end{array}
\right.
$$
\end{corollary}

From Lemma \ref{coclass}, one directly obtains the following.
\begin{proposition} A simply connected non-symmetric
$3$-dimensional homogeneous manifold is normal if and only if all
its first conjugate points of a fixed point are not isotropic.
\end{proposition}

\subsection{Tangent conjugate locus} Denote by $s_{p},$ $p\in {\mathbb N},$ the smooth map $s_{p}: ]0,
\pi[\to ]2p\pi, 2(p+1)\pi[$ where $s_{p}(\theta)\in {\mathcal
Z}^{+}(f_{\theta}),$ that is, $s_{p}(\theta)$ is the (unique)
solution of the equation $\tan \frac{s}{2} = \mu(\theta)s,$ such
that $s_{p}(\theta)\in ]2p\pi,2(p+1)\pi[.$ Then, as one can see in
{\it Fig 1} and {\it Fig 2}, we get:
\begin{enumerate}
\item[{\rm a)}] If
$\kappa<\tau^{2},$ $\lim_{\theta \to 0^{+}}s_{p}(\theta) =
\lim_{\theta \to \pi^{-}}s_{p}(\theta) = 2p\pi,$
$s_{p}(\frac{\pi}{2})$ is a maximum for $s_{p}$ on $]0,\pi[$ and
$s_{p}(]0,\pi[)\subset ]2p\pi,(2p+1)\pi[.$
\item[{\rm b)}] If $\kappa >
\tau^{2},$ $\lim_{\theta \to 0^{+}}s_{p}(\theta) = \lim_{\theta
\to \pi^{-}}s_{p}(\theta) = 2(p+1)\pi,$ $s_{p}(\frac{\pi}{2})$ is
a minimum and $s_{p}(]0,\pi[)\subset ](2p+1)\pi,2(p+1)\pi[.$
\end{enumerate}

We can extend $s_{p}$ to $[0,\pi]$ taking $s_{p}(0) = s_{p}(\pi) =
2p\pi,$ if $\kappa <\tau^{2},$ and $s_{p}(0) = s_{p}(\pi) =
2(p+1)\pi,$ if $\kappa >\tau^{2}.$ Moreover, if
$\tau^{2}-\kappa\to +\infty,$ then $s_{p}(\theta) \to (2p
+1)\pi^{-}$ and if $\tau^{2}-\kappa\to -\infty,$ then
$s_{p}(\theta) \to (2p + 1)\pi^{-},$ for all $\theta \in ]0,\pi[.$

Denote by ${\rm conj}_{\rm Isot}(M^{3}(\kappa,\tau))$ the
isotropic tangent conjugate locus of $M^{3}(\kappa,\tau).$ Then,
using Theorem \ref{pral1} and Proposition \ref{first}, we get

\begin{lemma}\label{lconj} The tangent conjugate locus ${\rm conj}(M^{3}(\kappa,
\tau))$ is the union of the following regular surfaces of
revolution:
$$
{\rm conj}(M^{3}(\kappa,\tau))= \left\{
\begin{array}{l}
 \bigcup_{p\in
\mathbb N}{\mathcal S}^{1}_{p}\cup {\mathcal S}^{2}_{p},\;\;
\mbox{if}\;\;
\kappa < \tau^{2},\\[0.6pc]
{\mathcal S}^{2}_{0} \cup \Big ( \bigcup_{p\in \mathbb N}{\mathcal
S}^{1}_{p}\cup {\mathcal S}^{2}_{p}\Big ),\;\; \mbox{if} \;\;
\kappa>\tau^{2},
\end{array}
\right.
$$
and
\[
{\rm conj}_{\rm Isot}(M^{3}(\kappa,\tau))  =   \bigcup_{p\in
{\mathbb N}}\Big ({\mathcal S}^{1}_{p} \setminus
\{(0,0,\pm\frac{2p\pi}{\tau})\}\Big ),
\]
where ${\mathcal S}^{1}_{p} =
\{\frac{2p\pi}{\sqrt{\lambda(\theta)}}u(\theta,\phi)\},$
${\mathcal S}^{2}_{0} =
\{\frac{{s}_{0}(\theta)}{\sqrt{\lambda(\theta)}}u(\theta,\phi)\}$
and ${\mathcal S}_{p}^{2}  =
\{\frac{{s}_{p}(\theta)}{\sqrt{\lambda(\theta)}}u(\theta,\phi)\},$
for each $p\in {\mathbb N},$ and for all $\theta\in
I_{M^{3}(\kappa,\tau)}$ and $\phi\in [0,2\pi].$
\end{lemma}

\begin{remark}{\rm ${\mathcal S}^{1}_{p}$ is the surface
generated by revolving the curve $\alpha^{1}_{p}(\theta) =
(2p\pi/\sqrt{\lambda(\theta)})u(\theta),$ for $p\in {\mathbb N}$
and $\theta\in I_{M^{3}(\kappa,\tau)},$ ${\mathcal S}^{2}_{p}$ by
revolving the curve $\alpha^{2}_{p}(\theta) =
(s_{p}(\theta)/\sqrt{\lambda(\theta)})u(\theta)$ and ${\mathcal
S}^{2}_{0},$ for $p= 0.$ All these surfaces are regular. In fact,
using that $\lim_{\theta\to 0^{+}}s_{p}'(\theta) = \lim_{\theta\to
\pi^{-}}s_{p}'(\theta)= 0$ and putting ${\mathfrak m} =
{\RR}^{3}[x,y,z],$ where $x,y,z$ are the cartesian coordinates
with respect to $\{e_{1},e_{2}, e_{3}\},$ we get
$$
\begin{array}{l}
\lim_{\theta\to 0^ {+}}(\alpha^{1}_{p})'(\theta) = -
\lim_{\theta\to\pi^{-}}(\alpha^{1}_{p})'(\theta) =
(\frac{2p\pi}{\tau}, 0 ,0)\\[0.6pc]
\lim_{\theta\to 0^ {+}}(\alpha^{2}_{p})'(\theta) = -
\lim_{\theta\to\pi^{-}}(\alpha^{2}_{p})'(\theta) =
(\frac{2p\pi}{\tau}, 0
,0),\;\;\mbox{if}\;\;\kappa<\tau^{2},\\[0.6pc]
\lim_{\theta\to 0^ {+}}(\alpha^{2}_{p})'(\theta) = -
\lim_{\theta\to\pi^{-}}(\alpha^{2}_{p})'(\theta) =
(\frac{2(p+1)\pi}{\tau}, 0 ,0), \;\; \mbox{if}\;\;\kappa>
\tau^{2}.
\end{array}
$$
}
\end{remark}

For the Heisenberg group $H_{3}(\tau),$ Lemma \ref{lconj} gives
\begin{proposition} We have:
\begin{enumerate}
\item[{\rm (i)}] The tangent conjugate locus ${\rm conj}(H_{3}(\tau))$ of
$H_{3}(\tau)$ is given by the union the planes $\Pi^{\pm}_{p}:z =
\pm\frac{2\pi p}{\tau},$ $p\in {\mathbb N},$ and the surfaces of
revolution parameterised as
\[
\vec{x}_{p}^{\pm}(\theta,\phi) = \frac{s_{p}(\theta)}{\tau} (\tan
\theta\cos \phi, \tan \theta \sin \phi,\pm 1),\;\;\;\; \theta\in
[0,\frac{\pi}{2}[.
\]
\item[{\rm (ii)}] $
{\rm conj}_{\rm Isot}(H_{3}(\tau))  =   \bigcup_{p\in {\mathbb
N}}\Big (\Pi_{p}^{\pm}\setminus \{(0,0,\pm\frac{2p\pi}{\tau})\}\big
).$
\item[{\rm (iii)}] The first tangent conjugate locus ${\rm conj}^{1}(H^{3}(\tau))$
of $H_{3}(\tau)$ are the planes $\Pi_{1}^{\pm}: z  =
\pm\frac{2\pi}{\tau}.$
\item[{\rm (iv)}] The first conjugate points to the origin are all
isotropic up to the points $\gamma(\pm \frac{2\pi}{\tau}),$ where
$\gamma$ is the Hopf fiber through the origin.
\end{enumerate}
\end{proposition}

From Theorem \ref{pral1}, $(x,y,z)$ belonging to ${\rm
conj}(M^{3}(\kappa,\tau))$ satisfies
\[
x=\frac{s\sin\theta\cos\phi}{\sqrt{\lambda(\theta)}},\;\;\;
y=\frac{s\sin\theta\sin\phi}{\sqrt{\lambda(\theta)}},\;\;\;
z=\frac{s\cos\theta}{\sqrt{\lambda(\theta)}},
\]
for $\theta \in I_{M^{3}(\kappa,\tau)}$ and $s\in{\mathcal
Z}^{+}(f_{\theta}).$ Then, we have $x^{2} + y^{2} + z^{2} =
\frac{s^{2}}{\lambda(\theta)}$ and for $\kappa\neq 0,$ one obtains
$\lambda(\theta) = \frac{\kappa s^{2}}{s^{2} -
(\tau^{2}-\kappa)z^{2}}.$ Hence, we get $\kappa(x^{2} + y^{2}) +
\tau^{2}z^{2} = s^{2}.$ It implies that the surfaces ${\mathcal
S}_{p}^{1},$ for $p\in {\mathbb N},$ in Lemma \ref{lconj} are
ellipsoids if $\kappa >0,$ and hyperboloids of two sheets if
$\kappa< 0,$ and ${\mathcal S}^{2}_{p}$ are surfaces of revolution
generated by revolving the curve $\kappa x^{2} + \tau^{2} z^{2} =
s^{2}_{p}(\theta)$ about the $z$-axis for $p\in {\mathbb N}\cup
\{0\}.$ Therefore, one can conclude with the following results:
\begin{proposition} We have:
\begin{enumerate}
\item[{\rm (i)}] The tangent conjugate locus ${\rm conj}(S^{3}(\kappa,\tau))$ of
$S^{3}(\kappa,\tau)$ is given by the union of:
\begin{enumerate}
\item[{\rm (a)}] the ellipsoids ${\mathcal E}_{p}: \;\kappa(x^{2} + y^{2}) + \tau^{2}z^{2} =
4p^{2}\pi^{2},$ $p\in {\mathbb N};$
\item[{\rm (b)}] the surfaces of revolution ${\mathcal S}_{p},$ for all $p\in {\mathbb N},$ generated by revolving the
curve $\kappa x^{2}(\theta) + \tau^{2} z^{2}(\theta) =
s^{2}_{p}(\theta)$ about the $z$-axis and moreover, for $p=0,$
when $\kappa>\tau^{2}.$
\end{enumerate}
\item[{\rm (ii)}] ${\rm conj}_{\rm Isot}(S^{3}(\kappa,\tau))=   \bigcup_{p\in {\mathbb N}}\Big
({\mathcal E}_{p}\setminus \{(0,0,\pm\frac{2p\pi}{\tau})\}\big ).$
\item[{\rm (iii)}] If $\kappa<\tau^{2},$ ${\rm conj}^{1}(S^{3}(\kappa,\tau)) ={\mathcal
E}_{1}.$
\item[{\rm (iv)}] If $\kappa>\tau,$ ${\rm conj}^{1}(S^{3}(\kappa,\tau))= {\mathcal
S}_{0}.$
\end{enumerate}
\end{proposition}

\begin{proposition} We have:
\begin{enumerate}
\item[{\rm (i)}] The tangent conjugate locus ${\rm conj}(\widetilde{SL_{2}}(\kappa,\tau))$ of
$\widetilde{SL_{2}}(\kappa,\tau)$ is given by the union of:
\begin{enumerate}
\item[{\rm (a)}] the hyperboloids of two sheets ${\mathcal H}_{p}: \;\kappa(x^{2} + y^{2}) + \tau^{2}z^{2} =
4p^{2}\pi^{2},$ $p\in {\mathbb N};$
\item[{\rm (b)}] the surfaces of revolution ${\mathcal S}_{p},$ for all $p\in {\mathbb N},$ generated by revolving the
curve $\kappa x^{2}(\theta) + \tau^{2} z^{2}(\theta) =
s^{2}_{p}(\theta)$ about the $z$-axis.
\end{enumerate}
\item[{\rm (ii)}] ${\rm conj}_{\rm Isot}(\widetilde{SL_{2}}(\kappa,\tau))=   \bigcup_{p\in {\mathbb N}}\Big
({\mathcal H}_{p}\setminus \{(0,0,\pm\frac{2p\pi}{\tau})\}\big ).$
\item[{\rm (iii)}] ${\rm conj}^{1}(\widetilde{SL_{2}}(\kappa,\tau)) ={\mathcal
H}_{1}.$
\item[{\rm (iv)}] The first conjugate points to the origin are all
isotropic up to the points $\gamma(\pm \frac{2\pi}{\tau}),$ where
$\gamma$ is the Hopf fiber through the origin.
\end{enumerate}
\end{proposition}

\end{document}